\documentclass[14pt]{article}
\usepackage{amssymb}
\usepackage{amsmath,amsfonts,amssymb,theorem,graphics,graphicx,epsfig,
    color,verbatim}
\usepackage{makeidx}

\makeindex

\newtheorem{Lem}{Lemma}

\newtheorem{Thm}{Theorem}

\newtheorem{Prop}{Proposition}

{\theorembodyfont{\normalfont}

\newtheorem{Rem}{Remark}

}

\def\BEN{\begin{enumerate}}
\def\BI{\begin{itemize}}
\def\EEN{\end{enumerate}}
\def\EI{\end{itemize}}

\def\beq{\begin{eqnarray}}
\def\eeq{\end{eqnarray}}

\def\beq{\begin{eqnarray}}
\def\eeq{\end{eqnarray}}
\def\bea{\begin{eqnarray*}}
\def\eea{\end{eqnarray*}}
\def\be{\begin{enumerate}}
\def\ee{\end{enumerate}}
\def\bi{\begin{itemize}}
\def\ei{\end{itemize}}

\input{tcilatex}

\begin{document}

\title{Limit theorems for additive functionals of stationary fields, under
integrability assumptions on the higher order spectral densities\thanks{
Partly supported by the Welsh Institute of Mathematics and Computational
Sciences and the Commission of the European Communities grant
PIRSES-GA-2008-230804 within the programme `Marie Curie Actions'. }}
\author{Florin Avram\thanks{%
Dept. de Math., Universite de Pau, France, E-mail: Florin.Avram@univ-Pau.fr}%
, Nikolai Leonenko\thanks{
Cardiff School of Mathematics, Cardiff University, Senghennydd Road,
Cardiff, CF24 4AG, UK, Email: LeonenkoN@Cardiff.ac.uk}, Ludmila Sakhno%
\thanks{%
Dept. of Probability, Statistics and Actuarial Mathematics, Kyiv National
Taras Shevchenko University, Ukraine, Email:lms@univ.kiev.ua}}
\maketitle

\begin{abstract}



\noindent We prove central limit theorems for additive functionals of
stationary fields under integrability conditions on the higher-order
spectral densities, which are derived using the H\"{o}%
lder-Young-Brascamp-Lieb inequality.

\ 


\end{abstract}

\section{Introduction}

\noindent{\textbf{Motivation.}} Consider a real measurable stationary in the
strict sense random field $X_{t} $, $t\in \mathbb{R}^{d},$ with ${\mathbb{E}}%
X_{t}=0,$ and ${\mathbb{E}}|X_{t}|^{k}<\infty$, $k=2,3,...$.

\noindent{\textbf{Assumption A:}} We will assume throughout the existence of
all order cumulants $c_{k}(t_{1},t_{2},...,t_{k})$ for our stationary random
field $X_{t}$, and also that they are representable as Fourier transforms of
``cumulant spectral densities'' \newline
$f_{k}(\lambda _{1},...,\lambda _{k-1})\in L_{1}(\mathbb{R}^{d(k-1)}),$ $%
k=2,3,...$, i.e: 
\begin{equation*}
c_{k}(t_{1},t_{2},...,t_{k})=c_{k}(t_{1}-t_{k},..,t_{k-1}-t_{k},0)=
\end{equation*}
\begin{equation*}
=\int_{\lambda _{1},...,\lambda _{k-1}\in \mathbb{R}^{d(k-1)}}e^{i%
\sum_{j=1}^{k-1}\lambda _{j}(t_{j}-t_{k})}f_{k}(\lambda _{1},...,\lambda
_{k-1})\,d\lambda _{1}...d\lambda _{k-1}.
\end{equation*}
\ \ \textbf{Note:} The functions $f_{k}(\lambda _{1},...,\lambda _{k-1})$
are symmetric and may be complex valued in general.

Central limit theorems for stationary fields have been derived traditionally
starting with the simplest cases of Gaussian or moving average processes,
via the method of moments based on explicit computations of the spectral
densities. We are able to treat here general stationary fields, by making
use of the powerful H\"{o}lder-Young-Brascamp-Lieb (HYBL) inequality.
Discussion of different approaches for derivation of CLT for stationary
processes and fields can be found, for example, in \cite{ALS}.

\noindent {\textbf{The problem:}} Let the random field $X_{t}$ be observed
over a sequence $K_{T}$ of increasing dilations of a bounded convex set $K$
of positive Lebesgue measure $|K|>0$, containing the origin, i.e. 
\begin{equation*}
K_{T}=TK,\quad T\rightarrow \infty .
\end{equation*}
Note that $|K_{T}|=T^{d}|K|$.

We investigate the asymptotic normality of the integrals 
\begin{equation}
S_{T}=\int_{t\in K_{T}}X_{t}dt  \label{int}
\end{equation}
and the integrals with a some weight function 
\begin{equation}
S_{T}^{w}=\int_{t\in K_{T}}w(t)X_{t}dt  \label{wint}
\end{equation}
as $T\rightarrow \infty $, without imposing any extra assumption on the
structure of the field such as linearity, etc. We will not also introduce
any kind of mixing conditions. We will establish central limit theorems for $%
S_{T}$ and $S_{T}^{w}$ , appropriately normalized, by the method of moments.
Namely, we will consider the cumulants of integrals (\ref{int}) and (\ref%
{wint}), represent them in the spectral domain, and evaluate their
asymptotic behavior basing on some analytic tools provided by harmonic
analysis. \ In such a way, via the spectral approach, all conditions needed
to prove the results will be concerned with integrability of the spectral
densities $f_{k}(\lambda _{1},...,\lambda _{k-1}),$ $k=2,3,...$.%

Taking consideration of the cumulants of $S_{T}$ (or $S_{T}^{w}$) in the
spectral domain one is lead to deal with some kind of convolutions of
spectral densities with particular kernel functions (see formulas for the
cumulants (\ref{cum}) and (\ref{cumw}) below). Similar convolutions have
been studied in the series of papers \ \ \cite{a} - \cite{ALS}, under the
name of Fejer matroid/graph integrals.

Estimates for this kind of convolutions follow from the H\"{o}%
lder-Young-Bras\-camp-Lieb inequality which, under prescribed conditions on
the integrability indices for a set of functions $f_{i}\in L_{p_{i}}(S,d\mu
) $, $i=1,...,n$, allows to write upper bounds for the integrals of the form 
\begin{equation}
\int_{S^{m}}\prod_{i=1}^{k}f_{i}(l_{i}(x_{1},...,x_{m}))\prod_{j=1}^{m}\mu
(dx_{j})  \label{holder}
\end{equation}%
with $l_{i}:S^{m}\rightarrow S$ being linear functionals (where $S$ may be
either torus $\left[ -\pi ,\pi \right] ^{d},$ $Z^{d},$ or $\mathbb{R}^{d}$
endowed with the corresponding Haar measure $\mu (dx)$).

An even more powerful tool, which we will need in this paper, is provided by
the nonhomogeneous H\"{o}lder-Young-Brascamp-Lieb inequality, which covers
the case when the above functions $f_{i}$ are defined over the spaces of
different dimensions: $f_{i}:S^{n_{i}}\rightarrow \mathbb{R}$ (see Appendix
A).

\noindent{\textbf{Contents:}} We state limit theorems for the integrals (\ref%
{int}) and (\ref{wint}) in Sections 2 and \ref{non-h} respectively, with
discussion of the assumptions used and of some possible applications. The
example of Gaussian fields is discussed in Section \ref{s:Ga}, and an
invariance principle provided in Section \ref{s:inv}. The H\"{o}%
lder-Young-Brascamp-Lieb inequality used to prove our results is presented
in Appendix A.

\section{Main results and discussion}

\label{s:intro}

Given a sequence $K_{T}$ of increasing dilations of a bounded convex set $K$
of positive Lebesgue measure $|K|>0$, containing the origin, let us consider
the uniform distribution on $K_{T}$ with the density 
\begin{equation*}
p_{K_{T}}(t)=\frac{1}{|K_{T}|}1_{\left\{ t\in K_{T}\right\} },\text{ }t\in 
\mathbb{R}^{d},
\end{equation*}
and characteristic function 
\begin{equation*}
\phi _{T}(\lambda )=\int_{\mathbb{R}^{d}}p_{K_{T}}(t)e^{it\lambda }dt=\frac{1%
}{|K_{T}|}\int_{K_{T}}e^{it\lambda }dt,\text{ }\lambda \in \mathbb{R}^{d}.
\end{equation*}

Define the Dirichlet type kernel 
\begin{equation}
\Delta _{T}(\lambda )=\int_{t\in K_{T}}e^{it\lambda }dt=|K_{T}|\phi
_{T}(\lambda ),\text{ }\lambda \in \mathbb{R}^{d}.  \label{ker}
\end{equation}

Denote 
\begin{equation}
\Delta _{1}(\lambda )=\int_{t\in K}e^{it\lambda }dt,\text{ }\lambda \in 
\mathbb{R}^{d}.
\end{equation}

We will need the following assumption:

\bigskip

\noindent {\textbf{Assumption K:}} The bounded convex set $K$ is such that: 
\begin{equation*}
C_{p}(K):=|\!|\Delta _{1}(\lambda )|\!|_{p}=\left( \int_{\mathbb{R}%
^{d}}|\Delta _{1}(\lambda )|^{p}d\lambda \right) ^{1/p}<\infty ,\quad
\forall p>p_{\ast }\geq 1.
\end{equation*}

\begin{Rem}
Assumption K and scaling imply 
\begin{equation}
|\!|\Delta _{T}(\lambda )|\!|_{p}=T^{d(1-1/p)}C_{p}(K).  \label{deltaTnorm}
\end{equation}
\end{Rem}

\begin{Rem}
The constants $C_{p}(K)$ and $p_{\ast }$ in Assumption K depend on Gaussian
curvature of the set $K$. This fact goes back to Van der Corput when $d=2$
-- see Herz (1962), Sadikova (1966), and Stein (1986) for extensions and
further references.
\end{Rem}

\bigskip

The explicit formula for $C_{p}(K)$ when $K$ is a cube: $K=[-1/2,1/2]^{d},$
is known: $C_{p}(K)=C_{p}^{d},$ where $C_{p}=\left( 2\int_{R}|\frac{\sin (z)%
}{z}|^{p}dz\right) ^{\frac{1}{p}},\quad \forall p>1.$ Note that in this case 
$p_{\ast }=1,$ and $C_{p_{1}}>C_{p_{2}}$ for $p_{1}<p_{2}$. For a ball $%
K_{T}=B_{T}=\{t\in R^{d}:\left\Vert t\right\Vert \leq T/2\}$ it is known
that 
\begin{equation*}
\Delta _{T}(\lambda )=\int_{B_{T}}e^{it\lambda }dt=\left( 2\pi \frac{T}{2}%
\right) ^{\frac{d}{2}}J_{d/2}\left( \left\Vert \lambda \right\Vert \frac{T}{2%
}\right) /\left\Vert \lambda \right\Vert ^{d/2},\quad \lambda \in \mathbb{R}%
^{d},
\end{equation*}%
where $J_{\nu }(z)$ is the Bessel function of the first kind and order $\nu
, $ and 
\begin{equation*}
C_{p}(K)=\left( 2\pi \right) ^{\frac{d}{2}}\;2^{^{-d(1-\frac{1}{p}%
)}}\left\vert s(1)\right\vert ^{1/p}\left( \int_{0}^{\infty }\rho
^{d-1}\left\vert \frac{J_{\frac{d}{2}}(\rho )}{\rho ^{d/2}}\right\vert
^{p}d\rho \right) ^{1/p},\ p>\frac{2d}{d+1},
\end{equation*}%
where $\left\vert s(1)\right\vert $ is the surface area of the unit ball in $%
\mathbb{R}^{d},~d\geq 2.$ In this case $p_{\ast }=\frac{2d}{d+1}>1,d\geq 2.$%
\bigskip

The derivation of the central limit theorem for the integrals (\ref{int})
will be based on the above estimates for the norms of functions $\Delta
_{T}(\lambda )$ and the important property of these functions stated in the
next lemma.

\begin{Lem}
The function 
\begin{equation*}
\Phi _{T}^{(2)}(\lambda )=\frac{1}{(2\pi )^{d}\left| K\right| T^{d}}\left|
\int_{t\in K_{T}}e^{it\lambda }dt\right| ^{2}=\frac{1}{(2\pi )^{d}\left|
K\right| T^{d}}\Delta _{T}(\lambda )\Delta _{T}(-\lambda ),\text{ }\lambda
\in \mathbb{R}^{d}
\end{equation*}
possesses the kernel properties (or is an approximate identity for
convolution): 
\begin{equation}
\int_{\mathbb{R}^{d}}\Phi _{T}^{(2)}(\lambda )d\lambda =1,  \label{ker1}
\end{equation}
and for any $\varepsilon >0$ when $T\rightarrow \infty $%
\begin{equation}
\lim \int_{\mathbb{R}^{d}\diagdown \varepsilon K}\Phi _{T}^{(2)}(\lambda
)d\lambda =0.  \label{ker2}
\end{equation}
\end{Lem}

\textit{Proof.}

The first relation (\ref{ker1}) follows from (\ref{ker}) and Plancherel
theorem. From Hertz(1962) and Sadikova(1966) one derives the following
assertion: if $\ K$ is a convex set and $\partial ^{(d-1)}\left\{ K\right\} $
is its surface area, then for any $\epsilon >0$%
\begin{equation*}
\int_{\left\| \lambda \right\| >\epsilon }\left| \int_{t\in K}e^{it\lambda
}dt\right| ^{2}d\lambda \leq \frac{8}{\epsilon }\partial ^{(d-1)}\left\{
K\right\} \left[ \int\limits_{0}^{\pi }\sin ^{d}zdz\right] ^{-1}
\end{equation*}
is valid. This inequality and homothety properties yields the second
relation (\ref{ker2}), see also Ivanov and Leonenko (1986), p.25).\bigskip

The cumulant of order $k\geq 2$ of the normalized integral $S_{T}$ is of the
form 
\begin{eqnarray*}
I_{T}^{(k)} &=&cum_{k}\left\{ \frac{S_{T}}{T^{d/2}},...,\frac{S_{T}}{T^{d/2}}%
\right\} \\
&=&\frac{1}{T^{dk/2}}\int_{t\in K_{T}}...\int_{t\in
K_{T}}c_{k}(t_{1}-t_{k},..,t_{k-1}-t_{k},0)dt_{1}...dt_{k} \\
&=&\frac{1}{T^{dk/2}}\int_{\lambda _{1},...,\lambda _{k-1}\in \mathbb{R}%
^{d(k-1)}}f_{k}(\lambda _{1},...,\lambda _{k-1})
\end{eqnarray*}
\begin{equation}
\times \Delta _{T}(\lambda _{1})...\Delta _{T}(\lambda _{k-1})\Delta
_{T}\left( -\sum_{i=1}^{k-1}\lambda _{i}\right) \,d\lambda _{1}...d\lambda
_{k-1},  \label{cum}
\end{equation}
where $\Delta _{T}(\lambda )$ is the Dirichlet type kernel (\ref{ker}).

\bigskip

To evaluate the second-order cumulant $I_{T}^{(2)}$ we will need one more
assumption.

\bigskip

\noindent \textbf{Assumption B:} The second-order spectral density $%
f_{2}(\lambda )$ is bounded and continuous and 
\begin{equation*}
f_{2}(0)=\frac{1}{(2\pi )^{d}}\int_{\mathbb{R}^{d}}\left( {\mathbb{E}}%
X_{t}X_{0}\right) dt\neq 0.
\end{equation*}

\bigskip

Under the assumption B we obtain from (\ref{cum}) and Lemma 1 as $%
T\rightarrow \infty $%
\begin{equation*}
cum_{2}\left\{ \frac{S_{T}}{T^{d/2}},\frac{S_{T}}{T^{d/2}}\right\}
=Var\left\{ \frac{S_{T}}{T^{d/2}}\right\} =(2\pi )^{d}\left| K\right| \int_{%
\mathbb{R}^{d}}\Phi _{T}^{(2)}(\lambda )f_{2}(\lambda )d\lambda \rightarrow
\end{equation*}
\begin{equation}
\rightarrow (2\pi )^{d}\left| K\right| f_{2}(0).  \label{var}
\end{equation}

To evaluate the integral (\ref{cum}) for $k\geq 3$ we apply the {H\"{o}%
lder-Young-Brascamp-Lieb inequality (see Theorem A1 in Appendix A). }

Comparing (\ref{cum}) and l.h.s. of (GH), we have in (\ref{cum}): $H=\mathbb{%
R}^{d(k-1)}$ and $k+1$ functions $g_{1}=g_{2}=...=g_{k}=\Delta _{T}$ on $%
\mathbb{R}^{d}$, $g_{k+1}=f_{k}$ on $\mathbb{R}^{d(k-1)}$; linear
transformations in our case are as follows: for $x=(x_{1},...,x_{k})\in 
\mathbb{R}^{d(k-1)}$ $l_{j}(x)=x_{j},$ $j=1,...,k-1,$ $l_{k}(x)=%
\sum_{j=1}^{k-1}x_{j},$ $l_{k+1}(x)=\mathrm{Id}$ (identity on $\mathbb{R}%
^{d(k-1)}$).

Suppose there exists $z=(z_{1},...,z_{k+1})\in \lbrack 0,1]^{k+1}$ such that
condition (C1) of {Theorem A1 is satisfied: } 
\begin{equation}
d(z_{1}+...+z_{k})+d(k-1)z_{k+1}=d(k-1),  \label{z}
\end{equation}
with 
\begin{equation*}
z_{1}=...=z_{k}=\frac{1}{p_{1}},\text{ }z_{k+1}=\frac{1}{p_{k+1}},
\end{equation*}
where $p_{1}$ falls in the range for which Assumption K holds, and $p_{k+1}$
is the integrability index of the spectral density $f_{k}$, that is, suppose 
$f_{k}(\lambda _{1},...,\lambda _{k-1})\in L_{p_{k+1}}(\mathbb{R}^{d(k-1)}).$

Let us check that condition (C2) will be satisfied as well with such a
choice of $z=(z_{1},...,z_{k+1}).$ For $\forall V\subset \mathbb{R}^{d(k-1)}$
we must have 
\begin{equation}
\dim V\leq \sum_{j=1}^{k+1}z_{j}\dim (l_{j}(V)),  \label{checkC2}
\end{equation}
the r.h.s. is equal to $z_{1}\sum_{j=1}^{k}\dim (l_{j}(V))+z_{k+1}\dim
(l_{k+1}(V))$, where, with the above choice of the linear transformations,
we have $\dim (l_{j}(V))=d,$ $j=1,...,k,$ $\dim (l_{k+1}(V))=\dim (V),$ that
is, \ (\ref{checkC2}) becomes 
\begin{equation*}
\dim V\leq z_{1}kd+z_{k+1}\dim (V),
\end{equation*}
or, taking into account that $(z_{1},...,z_{k+1})$ have chosen to satisfy (%
\ref{z}), 
\begin{equation*}
\dim V\leq \frac{z_{1}kd}{1-z_{k+1}}=\frac{d(k-1)\left( 1-z_{k+1}\right) }{%
1-z_{k+1}}=d(k-1),
\end{equation*}
which holds indeed for $\forall V\subset \mathbb{R}^{d(k-1)}.$

Then applying the {H\"{o}lder-Young-Brascamp-Lieb inequality (and taking
into account\ (\ref{deltaTnorm})) we have }for some $C>0$ 
\begin{equation}
\left| I_{T}^{(k)}\right| \leq CT^{kd(1-\frac{1}{p_{1}})-\frac{kd}{2}%
}C_{p_{1}}^{k}(K)||f_{k}||_{p_{k+1}}  \label{ineq}
\end{equation}
for $\forall p_{1}>p_{\ast }\geq 1$ and $\ p_{k+1}$ \ satisfuing (\ref{z}).

If $p_{\ast }<2,$ we can chose $p_{1}=2$ and come to the bound 
\begin{equation}
\left| I_{T}^{(k)}\right| \leq CC_{2}^{k}(K)||f_{k}||_{p_{k+1}},
\label{ineqO}
\end{equation}
for such a choice of $p_{1},$ the corresponding index $p_{k+1}$ we obtain
from (\ref{z}):

\begin{equation}
p_{k+1}=\frac{2(k-1)}{k-2},\text{ }k\geq 3.  \label{cond}
\end{equation}

However, we are able to prove that, in fact, $\ I_{T}^{(k)}\rightarrow 0$ as 
$T\rightarrow \infty $ (that is, bound in (\ref{ineqO}) can be strengthen to
the form $o(1)$), requiring still $p_{\ast }<2$ and $f_{k}(\lambda
_{1},...,\lambda _{k-1})\in L_{p_{k+1}}(\mathbb{R}^{d(k-1)})$ with the same $%
p_{k+1}$ given by (\ref{cond}).

Indeed, let us chose in (\ref{z}) $\tilde{p}_{1}=...=\tilde{p}_{k-2}=2$ \
(that is, $\tilde{z}_{1}=...=\tilde{z}_{k-2}=\frac{1}{2}$) and $\tilde{p}%
_{k-1}=\tilde{p}_{k}$ be close but less than $2$ ($\tilde{z}_{k-1}=\tilde{z}%
_{k}$ close but more than $\frac{1}{2}$).

Then the bound (\ref{ineq}) becomes 
\begin{eqnarray}
\left| I_{T}^{(k)}\right| &\leq &CT^{d(1-\frac{2}{\tilde{p}_{k}}%
)}C_{2}^{k-2}(K)C_{\tilde{p}_{k}}^{2}(K)||f_{k}||_{\tilde{p}_{k+1}}  \notag
\\
&=&CT^{-\varepsilon d}C_{2}^{k-2}(K)C_{\tilde{p}_{k}}^{2}(K)||f_{k}||_{%
\tilde{p}_{k+1}},  \label{ineq2}
\end{eqnarray}
where $\varepsilon =\frac{2}{\tilde{p}_{k}}-1>0$ and corresponding $\tilde{p}%
_{k+1}$, obtained from (\ref{z}), will be such that 
\begin{equation}
\tilde{p}_{k+1}>p_{k+1}=\frac{2(k-1)}{k-2},\text{ }k\geq 3  \label{cond1}
\end{equation}
(note that we do not need here the exact expressions for $\tilde{p}_{k}$ and 
$\tilde{p}_{k+1}$).

Therefore, for the functions $f_{k}\in L_{\tilde{p}_{k+1}}$ we have $%
I_{T}^{(k)}\rightarrow 0$ as $T\rightarrow \infty ,$ for $k\geq 3.$

Remembering that we are interested in evaluating (\ref{cum}) for the
functions $f_{k}$ which are in $L_{_{1}}$ (as being spectral densities), we
summarize the above reasonings as follows:

(i) for $f_{k}\in L_{_{1}}\cap L_{p_{k+1}}$ we have obtained the bound (\ref%
{ineqO});

(ii) for $f_{k}\in L_{_{1}}\cap L_{\tilde{p}_{k+1}}$ we have obtained the
convergence $I_{T}^{(k)}\rightarrow 0$ as $T\rightarrow \infty .$

It is left to note that

(iii) $L_{_{1}}\cap L_{\tilde{p}_{k+1}}$ is dense in $L_{_{1}}\cap
L_{p_{k+1}}$ (see (\ref{cond1}))

\noindent to conclude that the convergence $I_{T}^{(k)}\rightarrow 0$ as $%
T\rightarrow \infty $ holds for functions from $L_{_{1}}\cap L_{p_{k+1}}$ as
well.

Indeed, for $f_{k}\in L_{_{1}}\cap L_{p_{k+1}}$ and $g_{k}\in L_{_{1}}\cap
L_{\tilde{p}_{k+1}}$ we can write 
\begin{equation*}
\left| I_{T}^{(k)}(f_{k})\right| \leq \left| I_{T}^{(k)}(f_{k}-g_{k})\right|
+\left| I_{T}^{(k)}(g_{k})\right| ,
\end{equation*}
where the first term can be made arbitrary small with the choice of $g_{k}$
in view of (i) and (iii), and the second term tends to zero in view of (ii).

Thus, the following central limit theorem is proved by method of cumulants,
with conditions formulated in terms of spectral densities.

\begin{Thm}
Suppose that Assumptions A, K with $p_{\ast }<2$, and B hold, and for $k\geq
3$%
\begin{equation}
f_{k}(\lambda _{1},...,\lambda _{k-1})\in L_{p_{k}}(\mathbb{R}^{d(k-1)}),
\label{COND1}
\end{equation}
where $p_{k}=\frac{2(k-1)}{k-2}.$ Then, as $T\rightarrow \infty $ 
\begin{equation}
\frac{S_{T}}{T^{d/2}}\overset{D}{\rightarrow }N(0,\sigma ^{2}),  \label{CLT}
\end{equation}
where $\sigma ^{2}=(2\pi )^{d}\left| K\right| f_{2}(0).$
\end{Thm}

\begin{Rem}
For balls and cubes the condition $p_{\ast }<2$ holds$.$
\end{Rem}

\begin{Rem}
As a consequence of the above theorem we can state that the CLT (\ref{CLT})
holds under Assumptions A, K with $p_{\ast }<2$, and B, if the spectral
densities $f_{k}\in L_{4}(\mathbb{R}^{d(k-1)})$, $k\geq 3$. However, Theorem
1 provides more refined conditions, showing that for the central limit
theorem to hold the index of integrability of higher order spectral densities%
$f_{k}$ can become smaller and smaller, approaching to $2$ as $k$ grows.
\end{Rem}

The next remark is about a possible condition of the convex sets in form of
the kernel property.

\begin{Rem}
One can assume that the function 
\begin{equation*}
\Phi _{T}^{(k)}(\lambda _{1},...,\lambda _{k-1})=\frac{1}{(2\pi
)^{d(k-1)}\left| K\right| ^{k-1}T^{d}}\Delta _{T}(\lambda _{1})...\Delta
_{T}(\lambda _{k-1})\Delta _{T}\left( -\sum_{i=1}^{k-1}\lambda _{i}\right)
\end{equation*}
has the kernel property on \ $\mathbb{R}^{d(k-1)}$\ for $k\geq 2:$%
\begin{equation}
\int_{\mathbb{R}^{d(k-1)}}\Phi _{T}^{(k)}(\lambda _{1},...,\lambda
_{k-1})d\lambda _{1}...d\lambda _{k-1}=1,  \label{ker3}
\end{equation}
and for any $\varepsilon >0$ when $T\rightarrow \infty $ 
\begin{equation}
\lim \int_{\mathbb{R}^{d(k-1)}\diagdown \varepsilon K^{k-1}}\Phi
_{T}^{(k)}(\lambda _{1},...,\lambda _{k-1})d\lambda _{1}...d\lambda _{k-1}=0.
\label{ker4}
\end{equation}
Note that (\ref{ker3}), (\ref{ker4}) hold for the rectangle $K=\left[ -\frac{%
1}{2},\frac{1}{2}\right] ^{d}$ (see, for instance, Bentkus and Rutkauskas
(1973) or Avram, Leonenko and Sakhno (2010) and the references therein). If
the higher-order spectral densities $f_{k}(\lambda _{1},...,\lambda
_{k-1}),k\geq 2$ are continuous and bounded and $f_{k}(0,...,0)\neq 0,$ then 
\begin{equation*}
I_{T}^{(k)}=\frac{(2\pi )^{d}\left| K\right| ^{k-1}}{T^{d\left( \frac{k}{2}%
-1\right) }}\int_{\mathbb{R}^{d(k-1)}}\Phi _{T}^{(k)}(\lambda
_{1},...,\lambda _{k-1})f_{k}(\lambda _{1},...,\lambda _{k-1})d\lambda
_{1}...d\lambda _{k-1}\sim
\end{equation*}
\begin{equation*}
\sim \frac{(2\pi )^{d}\left| K\right| ^{k-1}}{T^{d\left( \frac{k}{2}%
-1\right) }}f_{k}(0,...,0),
\end{equation*}
as $T\rightarrow \infty ,$ thus tend to zero for $k\geq 3,$ and the central
limit theorem, Theorem 1, follows.
\end{Rem}

\section{Gaussian fields \label{s:Ga}}

Let us consider how the above method for deriving Theorem 1 can be used in
the situation when the field \ $X(t)$ is a nonlinear transformation of a
Gaussian field. Note that this kind of limit theorems, often called in the
literature Breuer-Major theorems, have been addressed by many authors.
Recently, powerful theory based on Malliavin calculus was exploited in the
series of papers by Nualart, Ortiz-Lattore, Nourdin, Peccati, Tudor and
others to develop CLTs in the framework of Wiener Chaos via remarkable
fourth moment approach (see, for example, \cite{NourdinPeccati2009}, \cite%
{NualartPeccati2005}\ \ and references therein). We show how CLT can be
stated quite straightforwardly with the use of the {H\"{o}%
lder-Young-Brascamp-Lieb inequality. }

For a stationary Gaussian filed $X(t),t\in \mathbb{R}^{d},$ consider the
field $Y(t)=G(X(t)),t\in \mathbb{R}^{d}.$ For a quite broad class of
functions $G$, evaluation of asymptotic behavior of the normalized integrals 
$S_{T}=\int_{t\in K_{T}}Y(t)dt$ reduces to consideration of the integrals $%
\int_{t\in K_{T}}H_{m}(X(t))dt,$ with a particular $m$, where $H_{m}(x)$ is
the Hermite polynomial, $m$ is Hermite rank of $G$ (see, i.e., Ivanov and
Leonenko (1986), p.55).

To demonstrate the approach based on the use of the {H\"{o}%
lder-Young-Brascamp-Lieb inequality, we consider here only the case of
integrals } 
\begin{equation}
S_{T}=S_{T}(H_{2}(X(t)))=\int_{t\in K_{T}}H_{2}(X(t))dt,  \label{*}
\end{equation}%
where $H_{2}(x)=x^{2}-1.$

Suppose that the centered Gaussian field $X(t),t\in \mathbb{R}^{d},$ has a
spectral density $f(\lambda ),\lambda \in \mathbb{R}^{d}.$ Then we can write
the following Wiener-It\^{o} integral representation: 
\begin{equation}
H_{2}(X(t))=\int_{\mathbb{R}^{2d}}e^{i(x,\lambda _{1}+\lambda _{2})}\sqrt{%
f(\lambda _{1})}\sqrt{f(\lambda _{2})}W(d\lambda _{1})W(d\lambda _{2}),
\label{**}
\end{equation}%
where $W(\cdot )$ is the Gaussian complex white noise measure (with
integration on the hyperplanes $\lambda _{i}=\pm \lambda _{j},i,j=1,2,i\neq j
$, being excluded).  Applying the formulas for the cumulants of multiple
stochastic Wiener-It\^{o} integrals, we have that the spectral density of
the second order of the field (\ref{**}) is given by 
\begin{equation*}
g_{2}(\lambda )=\int_{\mathbb{R}^{d}}f(\lambda )f(\lambda +\lambda
_{1})d\lambda _{1},
\end{equation*}%
which is well defined if $f(\lambda )\in L_{2}(\mathbb{R}^{d})$, and this
condition guarantees also that the Assumption B holds.

Next, the cumulants of the normalized integral (\ref{*}) can be written in
the form 
\begin{eqnarray*}
I_{T}^{(k)} &=&cum_{k}\left\{ \frac{S_{T}}{T^{d/2}},...,\frac{S_{T}}{T^{d/2}}%
\right\} \\
&=&\frac{1}{T^{dk/2}}\int_{\lambda _{1},...,\lambda _{k-1}\in \mathbb{R}%
^{d(k-1)}}\Delta _{T}(\lambda _{1})...\Delta _{T}(\lambda _{k-1})\Delta
_{T}\left( -\sum_{i=1}^{k-1}\lambda _{i}\right)
\end{eqnarray*}%
\begin{equation}
\times \int_{\mathbb{R}^{d}}f(\lambda )f(\lambda +\lambda _{1})...f(\lambda
+\lambda _{1}+...+\lambda _{k-1})d\lambda \,d\lambda _{1}...d\lambda _{k-1}.
\label{cumH2}
\end{equation}%
Now we can repeat the same reasonings as those for the proof of Theorem 1 to
conclude that $I_{T}^{(k)}\rightarrow 0$ as $T\rightarrow \infty ,$ for $%
k\geq 3$, under the condition $\ f(\lambda )\in L_{2}(\mathbb{R}^{d}).$

Indeed, formula (\ref{z}) relating the integrability indices $p$ for $\Delta
_{T}(\lambda )$ and $q$ for $f(\lambda )$ becomes in this case of the
following form: $dk\frac{1}{p}+dk\frac{1}{q}=dk$, or $\frac{1}{p}+\frac{1}{q}%
=1.$ We need already $f(\lambda )$ to be in $L_{2}(\mathbb{R}^{d})$ for a
proper behavior of the second order cumulant, therefore, choosing $q=2$, we
can take $p$ to be equal $2$ as soon as $p_{\ast }<2$ in the Assumption K.

Thus, we derived the known result (see, for example, \cite{IvanovLeonenko}):

\begin{Prop}
If a stationary Gaussian filed $X(t),t\in \mathbb{R}^{d},$ has the spectral
density $f(\lambda )\in L_{2}(\mathbb{R}^{d})$ and Assumptions K with $%
p_{\ast }<2$ holds, then, as $T\rightarrow \infty $ 
\begin{equation}
\frac{S_{T}(H_{2}(X(t)))}{T^{d/2}}\overset{D}{\rightarrow }N(0,\sigma ^{2}),
\label{CLTH2}
\end{equation}%
where 
\begin{equation}
\sigma ^{2}=(2\pi )^{d}\left\vert K\right\vert \int_{\mathbb{R}%
^{d}}f^{2}(\lambda )d\lambda .  \label{sigma2}
\end{equation}
\end{Prop}

As we can see, when taking into consideration the spectral domain, the
application of the {H\"{o}lder-Young-Brascamp-Lieb inequality allows to}
provide a very simple proof. Note also that this kind of technique has been
used for linear sequences (which generalize Gaussian fields) as well \cite%
{AF}.

Moreover, requiring more regularity on spectral density $f(\lambda )$, we
are able to evaluate the rate of convergence (\ref{CLTH2}) in the following
way.

Let us consider $\check{S}_{T}=\frac{S_{T}(H_{2}(X(t)))}{(2\pi
)^{d}\left\vert K\right\vert f_{2}(0)T^{d/2}}.$ We have for $f(\lambda )\in
L_{2}(\mathbb{R}^{d})$ the convergence as $T\rightarrow \infty $%
\begin{equation}
\check{S}_{T}\overset{D}{\rightarrow }N\backsim N(0,1).  \label{CLTN01}
\end{equation}%
We can state stronger version for this approximation, namely, that the
convergence (\ref{CLTN01}) takes place with respect to the Kolmogorov
distance: 
\begin{equation}
d_{Kol}(\check{S}_{T},N)=\sup_{z\in \mathbb{R}}|P(\check{S}%
_{T}<z)-P(N<z)|\rightarrow 0,  \label{dkol}
\end{equation}%
and also we can provide an upper bound for $d_{Kol}(\check{S}_{T},N).$ For
this we apply the results from \cite{NourdinPeccati2009}: since $\check{S}%
_{T}$ is representable as a double stochastic Wiener-It\^{o} integral we can
use the Proposition 3.8 of \cite{NourdinPeccati2009} which is concerned with
normal approximation in second Wiener Chaos and gives upper bounds for the
Kolmogorov distance solely in terms of the fourth and second cumulants. This
bound is of the form 
\begin{equation}
d_{Kol}(\check{S}_{T},N)\leq \sqrt{\frac{1}{6}cum_{4}(\check{S}%
_{T})+(cum_{2}(\check{S}_{T})-1)^{2}}.
\end{equation}%
So, we need only to control the fourth cumulant of $\check{S}_{T}$ and this
can be done with the use of the {H\"{o}lder-Young-Brascamp-Lieb inequality.
Due to this inequality, analogously to our previous derivations, for }$%
f(\lambda )\in L_{q}(\mathbb{R}^{d}),$ $q>2,$ and $\Delta _{T}(\lambda )\in
L_{p}(\mathbb{R}^{d}),$ with $\frac{1}{p}+\frac{1}{q}=1,$ we can write 
\begin{equation*}
\left\vert cum_{k}(\check{S}_{T})\right\vert \leq CT^{kd(1-\frac{1}{p})-%
\frac{kd}{2}}C_{p}^{k}(K)||f||_{q}^{k}=CT^{kd(\frac{1}{q}-\frac{1}{2}%
)}C_{p}^{k}(K)||f||_{q}^{k},
\end{equation*}%
therefore, 
\begin{equation*}
d_{Kol}(\check{S}_{T},N)\leq Const\text{ }T^{-\frac{q-2}{q}d},
\end{equation*}%
where the constant depends on $K$ and $f.$ Thus, the rate of convergence to
the normal law depends on the index of integralbility of $f(\lambda )$, in
particular, for $f(\lambda )\in L_{4}(\mathbb{R}^{d})$ we obtain 
\begin{equation*}
d_{Kol}(\check{S}_{T},N)\leq Const\text{ }\frac{1}{T^{d/2}}.
\end{equation*}%
The above technique can be also used for deriving CLT for $%
S_{T}(H_{m}(X(t))) $ with $m>2.$ 


\section{An invariance principle \label{s:inv}}

Let us return now to the case of a general random field $X(t)$ of Assumption
A. In order to discuss the invariance principle for the situation above we
consider the multiparameter Brownian motion of Chentsov's type (see
Samorodnitsky and Taqqu (2004) for example), that is the zero mean Gaussian
random field $b(t),t\in \mathbb{R}^{d},$ such that

(i) $b(t)=0,$ if $t_{j}=0$ for at least one $j\in \left\{ 1,...,d\right\} ;$

(ii) ${\mathbb{E}}b(t_{1})b(t_{2})=\prod\limits_{j=1}^{d}\min \left\{
t_{1}^{(j)},t_{2}^{(j)}\right\} ,t_{l}=\left( t_{l}^{(j)},j=1,...,d\right)
,l\in \left\{ 1,2\right\} .$

We introduce the Gaussian process 
\begin{equation}
L_{K}(u)=\left( (2\pi )^{d}f_{2}(0)\right) ^{1/2}\int_{t\in u^{1/d}K}db(t),%
\text{ \ \ }u\in \left[ 0,1\right] ,  \label{inv1}
\end{equation}
with zero mean and covariance function 
\begin{equation*}
{\mathbb{E}}L_{K}(u_{1})L_{K}(u_{2})=(2\pi )^{d}f_{2}(0)\left| u_{1}^{\frac{1%
}{d}}K\cap u_{2}^{\frac{1}{d}}K\right| ,\text{ \ \ }u_{1},u_{1}\in \left[ 0,1%
\right] .
\end{equation*}

Note that for the ball $K=B_{1}=\{t\in \mathbb{R}^{d}:\left\| t\right\| \leq
1/2\}$%
\begin{equation*}
{\mathbb{E}}L_{B_{1}}(u_{1})L_{B_{1}}(u_{2})=(2\pi )^{d}f_{2}(0)\left|
B_{1}\right| \min \left\{ u_{1},u_{2}\right\} ,\text{ \ \ }u_{1},u_{2}\in %
\left[ 0,1\right] ,
\end{equation*}
where $\left| B_{1}\right| $ is the volume of the ball $B_{1}.$

If we assume that the stochastic process (\ref{inv1}) induces the
probabilistic measure $\mathit{P}$ in the space $C[0,1]$ of continuous
functions with the uniform topology, then one can prove the invariance
principle for the measures $\mathit{P}_{T},$ induced in the space $C[0,1]$
by the stochastic processes 
\begin{equation}
Y_{T}(u)=\frac{1}{T^{d/2}}\int_{t\in u^{1/d}K_{T}\,}X_{t}dt,\text{ \ }u\in %
\left[ 0,1\right] ,  \label{YT}
\end{equation}
that is the under conditions of Theorem 1 the measures $\mathit{P}_{T}$
converge weakly ($\Longrightarrow )$ to the Gaussian measure $\mathit{P}$ in
the space $C[0,1]$ \ as $T\rightarrow \infty $ (see Billingsley (1968) for
necessary definitions related to the convergence of probability measures).
This can be proved if we introduce the following assumption.

\bigskip

\noindent \textbf{Assumption }$\mathbf{\Phi }$\textbf{: }The function 
\begin{equation}
\Phi _{T_{1},T_{2}}^{(2)}(\lambda )=\frac{1}{(2\pi )^{d}\left| K\right|
(T_{2}^{d}-T_{1}^{d})}\left| \int_{t\in K_{T_{2}}\diagdown
K_{T_{1}}}e^{it\lambda }dt\right| ^{2},\text{ }\lambda \in \mathbb{R}^{d},
\label{Phi}
\end{equation}
has the kernel properties similar to (\ref{ker1}), (\ref{ker2}) for $%
T_{1}<T_{2},T_{1}\rightarrow \infty .$

\bigskip

Really, in this case one can check that the Kolmogorov's criterion: 
\begin{equation}
{\mathbb{E}}\left| Y_{T}(u_{2})-Y_{T}(u_{1})\right| ^{4}\leq const\left|
u_{2}-u_{1}\right| ^{2},\text{ }0\leq u_{1}\leq u_{2}\leq 1,  \label{moment4}
\end{equation}
of weakly compactness of probability measures $\left\{ \mathit{P}%
_{T}\right\} $ is satisfied (see again Billingsley (1968)).

Consider 
\begin{eqnarray}
{\mathbb{E}}\left| Y_{T}(v)-Y_{T}(u)\right| ^{4} &=&\frac{1}{T^{2d}}{\mathbb{%
E}}\left[ \int_{t\in v^{1/d}K_{T}\,\backslash u^{1/d}K_{T}}X_{t}dt\right]
^{4}  \notag \\
&=&\frac{1}{T^{2d}}\int_{\widetilde{K}_{T}^{4}\text{ }}{\mathbb{E}}\left[
X_{t_{1}}X_{t_{2}}X_{t_{3}}X_{t_{4}}\right] dt_{1}dt_{2}dt_{3}dt_{4}  \notag
\\
&=&\frac{1}{T^{2d}}\int_{\widetilde{K}_{T}^{4}\text{ }%
}[c_{4}(t_{1}-t_{4},t_{2}-t_{4},t_{3}-t_{4})  \notag \\
&&+c_{2}(t_{1}-t_{2})c_{2}(t_{3}-t_{4})+c_{2}(t_{1}-t_{3})c_{2}(t_{2}-t_{4})
\notag \\
&&+c_{2}(t_{1}-t_{4})c_{2}(t_{2}-t_{3})]dt_{1}dt_{2}dt_{3}dt_{4}  \notag \\
&=&I_{1}+I_{2}+I_{3}+I_{4}.  \label{sumi}
\end{eqnarray}
(We have denoted here $\widetilde{K}_{T}=v^{1/d}K_{T}\,\backslash
u^{1/d}K_{T}$.)

We can write 
\begin{eqnarray}
I_{1} &=&\frac{1}{T^{2d}}\int_{R^{3d}}f_{4}(\lambda _{1},\lambda
_{2},\lambda _{3})\prod_{j=1}^{3}\left[ \int_{\widetilde{K}_{T}\text{ }%
}e^{it_{j}\lambda _{j}}dt_{j}\right] \int_{\widetilde{K}_{T}\text{ }%
}e^{-it_{4}\sum_{j=1}^{3}\lambda _{j}}dt_{4}d\lambda _{1}d\lambda
_{2}d\lambda _{3}  \notag \\
&=&\frac{1}{T^{2d}}\int_{R^{3d}}f_{4}(\lambda _{1},\lambda _{2},\lambda
_{3})\prod_{j=1}^{3}\left[ \Delta _{\widetilde{K}_{T}}(\lambda _{j})\right]
\Delta _{\widetilde{K}_{T}}(\sum_{j=1}^{3}\lambda _{j})d\lambda _{1}d\lambda
_{2}d\lambda _{3}.  \label{I1}
\end{eqnarray}
Supposing $f_{4}\in L_{q}$, $\Delta _{\widetilde{K}_{T}}\in L_{q}$ for $p,q:4%
\frac{1}{p}+3\frac{1}{q}=3$ and applying the {H\"{o}lder-Young-Brascamp-Lieb
inequality we obtain } 
\begin{equation}
I_{1}\leq \frac{1}{T^{2d}}\left\| f_{4}\right\| _{q}\left\{ \left\| \Delta _{%
\widetilde{K}_{T}}\right\| _{p}\right\} ^{4}.  \label{I1<=}
\end{equation}
Choosing $p=2$ we get 
\begin{eqnarray*}
\left\{ \left\| \Delta _{\widetilde{K}_{T}}\right\| _{2}\right\} ^{4}
&=&\left\{ \left\{ \int_{R^{d}}\left| \int_{\widetilde{K}_{T}\text{ }%
}e^{it\lambda }dt\right| ^{2}d\lambda \right\} ^{1/2}\right\} ^{4}=\left\{
\int_{R^{d}}\left| \int_{\widetilde{K}_{T}\text{ }}e^{it\lambda }dt\right|
^{2}d\lambda \right\} ^{2} \\
&=&\left\{ (2\pi )^{d}|K|((Tv^{1/d})^{d}-(Tu^{1/d})^{d})\int_{R^{d}}\Phi
_{Tu^{1/d},Tv^{1/d}}^{(2)}(\lambda )d\lambda \right\} ^{2}.
\end{eqnarray*}
Therefore, under the assumption $f_{4}\in L_{3}(\mathbb{R}^{3d})$ (which is
covered by the assumptions of Theorem 1)

\begin{equation*}
I_{1}\leq const\text{ }(v-u)^{2}.
\end{equation*}

Next, consider 
\begin{equation*}
\int_{\widetilde{K}_{T}^{2}\text{ }}c_{2}(t_{1}-t_{2})dt_{1}dt_{2}=%
\int_{R^{d}}f_{2}(\lambda )\int_{\widetilde{K}_{T}^{2}\text{ }}e^{i\left(
t_{1}-t_{2}\right) \lambda }dt_{1}dt_{2}d\lambda =\int_{R^{d}}f_{2}(\lambda
)\left| \Delta _{\widetilde{K}_{T}}(\lambda )\right| ^{2}d\lambda .
\end{equation*}
Supposing $f_{2}(\lambda )$ to be bounded we get 
\begin{equation}
\int_{\widetilde{K}_{T}^{2}\text{ }}c_{2}(t_{1}-t_{2})dt_{1}dt_{2}\leq const%
\text{ }\int_{R^{d}}\left| \Delta _{\widetilde{K}_{T}}(\lambda )\right|
^{2}d\lambda =const\text{ }(2\pi )^{d}|K|T^{d}(v-u)  \label{Ij<=}
\end{equation}
which implies that each term $I_{j},$ $j=2,3,4$ in (\ref{sumi}) is bounded
by 
\begin{equation*}
const(2\pi )^{d}|K|(v-u)^{2}.
\end{equation*}

Hence, \ (\ref{moment4}) holds if we suppose that the second order spectral
density $f_{2}$ is bounded, $f_{4}\in L_{3}$ and $\Phi
_{T_{1},T_{2}}^{(2)}(\lambda )$ given by (\ref{Phi}) has the kernel
properties.

If the homogeneous random field $X_{t},t\in \mathbb{R}^{d},$ is second-order
isotropic (it means that the covariance function ${\mathbb{E}}%
X_{t}X_{s}=B\left( \left\| t-s\right\| \right) $ depends on the Euclidean
distance $\left\| t-s\right\| ,t,s\in \mathbb{R}^{d},$ and $%
K_{T}=B_{T}=\{t\in \mathbb{R}^{d}:\left\| t\right\| \leq T/2\}$ are balls,
then the condition (\ref{moment4}) and Assumption $\Phi $ concerning the
kernel properties of $\Phi _{T_{1},T_{2}}^{(2)}(\lambda )$ are satisfied. It
follows from the results by Leonenko and Yadrenko (1979) (see also Ivanov \
and Leonenko (1989), chapter 2), since for balls 
\begin{equation*}
{\mathbb{E}}\left[ \int_{T_{1}\leq \left\| t\right\| \leq T_{2}}X_{t}dt%
\right] ^{2}=\frac{4\pi ^{d}}{d\Gamma ^{2}(\frac{d}{2})}\gamma \left(
T_{2}^{d}-T_{1}^{d}\right) (1+o(1)),\text{ \ }T_{1}\rightarrow \infty ,
\end{equation*}
if 
\begin{equation*}
\int\limits_{0}^{\infty }z^{d-1}\left| B(z)\right| dz<\infty ,\text{ \ \ }%
\gamma =\int\limits_{0}^{\infty }z^{d-1}B(z)dz\neq 0.
\end{equation*}

We can summarize the above arguments in the next theorem.

\begin{Thm}
Suppose that Assumptions A, K, B and $\Phi $ hold, and $f_{4}(\lambda
_{1},\lambda _{2},\lambda _{3})\in L_{3}(\mathbb{R}^{3d}).$ Then the familly
of measures $\mathit{P}_{T},$ induced by the stochastic processes (\ref{YT})
is weakly compact in the space $C[0,1].$
\end{Thm}

Compiling now Theorem 1 and Theorem 2 we come to the following result.

\begin{Thm}
Suppose that conditions of Theorem 1 and, in addition, Assumptions $\Phi $
hold. Then \ $\mathit{P}_{T}\Longrightarrow \mathit{P}$ in $C[0,1]$, where
the measures $\mathit{P}_{T}$ and $\mathit{P}$ are induced by the stochastic
processes (\ref{YT}) and (\ref{inv1}) respectively.
\end{Thm}

\section{Non-homogeneous random fields \label{non-h}}

We discuss now the central limit theorem for non-homogeneous random fields
of special form.

\bigskip

\noindent {\textbf{Assumption C: }}Assume that a real (weight) function $%
w(t),t\in \mathbb{R}^{d},$ is (positively) homogeneous of degree $\beta ,$
that is for any $a>0$ there exists $\beta \in \mathbb{R}$, such that 
\begin{equation*}
w(at)=w(at_{1},...,at_{d})=a^{\beta }w(t),\text{ \ }t\in \mathbb{R}^{d}.
\end{equation*}

\ 

\noindent{\textbf{Assumption D: }}Assume that there exists 
\begin{equation*}
w_{1}(\lambda )=\int_{t\in K}w(t)e^{it\lambda }dt,\text{ }\lambda \in 
\mathbb{R}^{d}.
\end{equation*}

Under Assumptions C and D 
\begin{equation*}
w_{T}(\lambda )=\int_{t\in K_{T}}w(t)e^{it\lambda }dt=T^{d+\beta
}w_{1}(\lambda T),\text{ }\lambda \in \mathbb{R}^{d}.
\end{equation*}

\bigskip

\noindent \emph{Example 1}. The function $w_{1}(t)=\left\| t\right\| ^{\nu
},\nu \geq 0,$ is homogeneous of degree $\beta =\nu ,$ if $\nu >0.$ For
example if $d=1$, $K=[0,1]$ and $\nu \geq 0$ is an integer,we obtain 
\begin{equation*}
w_{1}(\lambda )=\int_{t\in \lbrack 0,1]}t^{\nu }e^{it\lambda }dt=\frac{1}{%
i^{\nu }}\frac{\partial ^{\nu }}{\partial \lambda ^{\nu }}\int_{t\in \lbrack
0,1]}e^{it\lambda }dt=\frac{1}{i^{\nu }}\frac{\partial ^{\nu }}{\partial
\lambda ^{\nu }}\frac{e^{i\lambda }-1}{i\lambda },\text{ }\lambda \in 
\mathbb{R}^{1}.
\end{equation*}

\bigskip

\noindent \emph{Example 2}. Another example of the homogeneous function of
degree $\beta >0,$ is $w_{2}(t)=\left| t_{1}+...+t_{d}\right| ^{\nu },$
where again $\beta =\nu ,$ if $\nu >0.$

\bigskip

\noindent \emph{Example 3}. The function $w_{3}(t)=\left| \left|
t_{1}\right| ^{\gamma }+..+.\left| t_{d}\right| ^{\gamma }\right| ^{\nu }$
is homogeneous of degree $\beta =\nu \gamma ,$ if $\nu >0,\gamma >0.$

\bigskip

\noindent \emph{Example 4}. All arithmetic, geometric and harmonic\ averages
of $\left| t_{1}\right| ,...,\left| t_{d}\right| $ are homogeneous functions
of degree one.

\bigskip

Under Assumption C we investigate below the asymptotic normality of
integrals 
\begin{equation*}
S_{T}^{w}=\int_{t\in K_{T}}w(t)X_{t}dt
\end{equation*}
as $T\rightarrow \infty $.

We denote 
\begin{equation}
W^{2}(T)=\int_{t\in K_{T}}w^{2}(t)dt=\frac{1}{(2\pi )^{d}}\int_{\mathbb{R}%
^{d}}\left| w_{T}(\lambda )\right| ^{2}d\lambda .  \label{W2Tdef}
\end{equation}
{\textbf{Assumption E: }}Let the finite measures 
\begin{equation*}
\mu _{T}(d\lambda )=\frac{\left| w_{T}(\lambda )\right| ^{2}d\lambda }{\int_{%
\mathbb{R}^{d}}\left| w_{T}(\lambda )\right| ^{2}d\lambda },\lambda \in 
\mathbb{R}^{d}
\end{equation*}
converge weakly to some finite measure $\mu (d\lambda ),$ and the spectral
density $f_{2}(\lambda )$ is positive on set $B$ $\subseteq \mathbb{R}^{d}$
of positive $\mu $-measure ($\mu (B)>0$).

We recall that the weak convergence of probability measures means that for
any continuous and bounded function $f(\lambda )$ as $T\rightarrow \infty $%
\begin{equation*}
\lim \int_{\mathbb{R}^{d}}f(\lambda )\mu _{T}(d\lambda )=\int_{\mathbb{R}%
^{d}}f(\lambda )\mu (d\lambda ).
\end{equation*}
Then we have that the variance 
\begin{equation*}
{\mathbb{E}}\left[ \frac{S_{T}^{w}}{W(T)}\right] ^{2}=\frac{1}{W^{2}(T)}%
\int_{\mathbb{R}^{d}}f(\lambda )\left[ \int_{t_{1}\in
K_{T}}w(t_{1})e^{it_{1}\lambda }dt_{1}\right] \overline{\left[
\int_{t_{2}\in K_{T}}w(t_{2})e^{it_{2}\lambda }dt_{2}\right] }d\lambda =
\end{equation*}
\begin{equation*}
=(2\pi )^{d}\int_{\mathbb{R}^{d}}f_{2}(\lambda )\mu _{T}(d\lambda
)\rightarrow (2\pi )^{d}\int_{\mathbb{R}^{d}}f_{2}(\lambda )\mu (d\lambda
)=\sigma ^{2}>0,
\end{equation*}
as $T\rightarrow \infty .$

It turns out that we need the following

{\textbf{Assumption F: }} 
\begin{equation*}
C_{p,w}(K):=||w_{1}(\lambda )||_{p}=\left( \int_{\mathbb{R}^{d}}\left|
\int_{t\in K}w(t)e^{it\lambda }dt\right| ^{p}d\lambda \right) ^{1/p}<\infty
,\quad \forall p>p_{\ast }\geq 1.
\end{equation*}

Then by scaling property we obtain the following formula: 
\begin{equation*}
|\!|w_{T}(\lambda )|\!|_{p}=T^{d(1-\frac{1}{p})+\beta }C_{p,w}(K),
\end{equation*}
and in particular 
\begin{equation}
W^{2}(T)=\int_{t\in K_{T}}w^{2}(t)dt=\frac{1}{(2\pi )^{d}}\int_{\mathbb{R}%
^{d}}\left| w_{T}(\lambda )\right| ^{2}d\lambda =\left[ \frac{1}{(2\pi )^{d}}%
T^{\frac{d}{2}+\beta }C_{2,w}(K)\right] ^{2}.  \label{W2T}
\end{equation}

Similar to the proof of Theorem 1 we obtain that the cumulant of order $%
k\geq 3$ is of the form 
\begin{equation*}
I_{T}^{(k)}=cum_{k}\left\{ \frac{S_{T}^{w}}{W(T)},...,\frac{S_{T}^{w}}{W(T)}%
\right\} =
\end{equation*}
\begin{equation*}
=\frac{1}{W(T)^{k}}\int_{t\in K_{T}}...\int_{t\in
K_{T}}w(t_{1})...w(t_{k})c_{k}(t_{1}-t_{k},..,t_{k-1}-t_{k},0)dt_{1}...dt_{k}=
\end{equation*}
\begin{equation*}
=\frac{1}{W(T)^{k}}\int_{\lambda _{1},...,\lambda _{k-1}\in \mathbb{R}%
^{d(k-1)}}w_{T}(\lambda _{1})w_{T}(\lambda _{2})...w_{T}(\lambda
_{k-1})w_{T}(-\lambda _{1}-...-\lambda _{k-1})\times
\end{equation*}

\begin{equation}
\times f_{k}(\lambda _{1},...,\lambda _{k-1})\,d\lambda _{1}...d\lambda
_{k-1},  \label{cumw}
\end{equation}
and then applying the {H\"{o}lder-Young-Brascamp-Lieb inequality with the
same reasonings as those used for derivation of the formula (\ref{ineq}) we
obtain }for some $C>0$ the bound 
\begin{equation*}
\left| I_{T}^{(k)}\right| \leq C\frac{T^{kd(1-\frac{1}{p_{1}})}}{T^{\frac{kd%
}{2}}C_{2,w}^{k}(K)}C_{p_{1},w}^{k}(K)|\!|f_{k}|\!|_{p_{k+1}}=CT^{-\nu
}C_{p_{1},w}^{k}(K)C_{2,w}^{-k}(K)|\!|f_{k}|\!|_{p_{k+1}},
\end{equation*}
where 
\begin{equation*}
\nu =kd\left( \frac{1}{2}-\left( 1-\frac{1}{p_{1}}\right) \right) .
\end{equation*}

Similar to the proof of the Theorem 1, from the condition $\nu >0$ we come
to the restrictions on $p_{1}$ and $p_{k+1}$, and, therefore, derive the
following

\begin{Thm}
If Assumptions A, C, D, E and F hold, and for $k\geq 3$%
\begin{equation*}
f_{k}(\lambda _{1},...,\lambda _{k-1})\in L_{p_{k}}(\mathbb{R}^{d(k-1)}),
\end{equation*}
where $p_{k}=\frac{2(k-1)}{k-2}.$ Then, as $T\rightarrow \infty $ 
\begin{equation*}
\frac{S_{T}^{w}}{W(T)}\overset{D}{\rightarrow }N(0,\sigma ^{2}),
\end{equation*}
where $\sigma ^{2}=(2\pi )^{d}\int_{\mathbb{R}^{d}}f_{2}(\lambda )\mu
(d\lambda ),$ and the finite measure $\mu $ is defined in assumption E.
\end{Thm}

This theorem can be applied to the statistical problem of estimation of
unknown coefficient of linear regression observed on the increasing convex
sets.

\bigskip

Analogously to Section 2, the invariance principle for the above situation
can be considered and Theorem 4 can be extended to the analog of Theorem 3.
We just point out the key steps here.

First, we note that for the monotonically increasing function $%
V(T):=W^{2}(T) $ (with $W^{2}(T)$ given by (\ref{W2Tdef})) there exists the
unique inverse function which we will denote $V^{(-1)}(T).$ Then we make the
modifications in the definitions of the processes (\ref{inv1}) and (\ref{YT}%
). The Gaussian process (\ref{inv1}) is defined now as the process 
\begin{equation*}
L_{K}^{w}(u)=\left( (2\pi )^{d}\int_{\mathbb{R}^{d}}f(\lambda )\mu (d\lambda
)\right) ^{1/2}\int_{t\in V^{(-1)}(u)K}db(t),\text{ \ \ }u\in \left[ 0,1%
\right] ,
\end{equation*}
with zero mean and the covariance function 
\begin{equation*}
{\mathbb{E}}L_{K}^{w}(u_{1})L_{K}^{w}(u_{2})=(2\pi )^{d}\int_{\mathbb{R}%
^{d}}f(\lambda )\mu (d\lambda )\left| V^{(-1)}(u_{1})K\cap
V^{(-1)}(u_{2})K\right| ,\,u_{1},u_{1}\in \left[ 0,1\right] .
\end{equation*}
Instead of (\ref{YT}) we consider the process 
\begin{equation}
Y_{T}^{w}(u)=\frac{1}{V(T)^{1/2}}\int_{t\in V^{(-1)}(u)K_{T}\,}w(t)X_{t}dt,%
\text{ \ }u\in \left[ 0,1\right] .  \label{YTw}
\end{equation}
Basing the proof of weak compactness of measures $\mathit{P}_{T}$ induced by
the stochastic processes (\ref{YTw}) on Kolmogorov's criterion (\ref{moment4}%
), we must check now that 
\begin{equation}
\frac{1}{V(T)^{2}}{\mathbb{E}}\left[ \int_{t\in V^{(-1)}(v)K_{T}\,\backslash
V^{(-1)}(u)K_{T}}w(t)X_{t}dt\right] ^{4}\leq const\left| v-u\right| ^{2},%
\text{ }0\leq u\leq v\leq 1.  \label{moment4w}
\end{equation}
The same derivations as those in Section 3 will lead to the expression for
the right hand side of (\ref{moment4w}) in the form of the sum $%
I_{1}+I_{2}+I_{3}+I_{4},$ where now the function $w(t)$ will be involved and
correspondingly in the formulas (\ref{I1}), (\ref{I1<=}) and (\ref{Ij<=}) $%
\Delta _{\widetilde{K}_{T}}(\lambda )$ will be changed for $\Delta _{%
\widetilde{K}_{T}}^{w}(\lambda )=\int_{t\in \widetilde{K}_{T}}w^{2}(t)dt$
with $\widetilde{K}_{T}$ being now of the form $\widetilde{K}%
_{T}=V^{(-1)}(v)K_{T}\,\backslash V^{(-1)}(u)K_{T}.$

Therefore, supposing $f_{2}$ to be bounded and $f_{4}\in L_{3},$ we come to
the following bound 
\begin{eqnarray}
&&{\mathbb{E}}\left[ \int_{t\in V^{(-1)}(v)K_{T}\,\backslash
V^{(-1)}(u)K_{T}}w(t)X_{t}dt\right] ^{4}  \notag \\
&\leq &const\left\{ \int_{R^{d}}\left| \Delta _{\widetilde{K}%
_{T}}^{w}(\lambda )\right| ^{2}d\lambda \right\} ^{2}  \notag \\
&=&const\left\{ \int_{R^{d}}\left| \int_{t\in \widetilde{K}%
_{T}}w(t)e^{it\lambda }dt\right| ^{2}d\lambda \right\} ^{2}  \notag \\
&=&const(2\pi )^{2d}\left\{ \int_{t\in \widetilde{K}_{T}}w^{2}(t)dt\right\}
^{2}.  \label{moment4w<=}
\end{eqnarray}
(Note that (\ref{moment4w<=}) can be compared with the formula (1.8.11) in 
\cite{IvanovLeonenko}, which gives more general result, namely, the bounds
for odd order higher moments).

Using (\ref{W2T}) we can derive 
\begin{eqnarray}
\int_{t\in \widetilde{K}_{T}}w^{2}(t)dt &=&\int_{t\in
V^{(-1)}(v)K_{T}\,\backslash V^{(-1)}(u)K_{T}}w^{2}(t)dt  \notag \\
&=&V(TV^{(-1)}(v))-V(TV^{(-1)}(u))  \notag \\
&=&T^{d+2\beta }(V(V^{(-1)}(v))-V(V^{(-1)}(u)))  \notag \\
&=&T^{d+2\beta }(v-u).  \label{intw2}
\end{eqnarray}
From (\ref{W2T}) we know also that $V(T)=(2\pi )^{-2d}T^{d+2\beta }const,$
which combined with (\ref{intw2}) and (\ref{moment4w<=}) gives (\ref%
{moment4w}). Therefore, weak compactness of measures $\mathit{P}_{T}$
induced by the stochastic processes (\ref{YTw}) takes place under the
conditions that the second order spectral density $f_{2}$ is bounded and the
fourth order spectral density $f_{4}$ is in $L_{3}.$

\bigskip

\section*{Appendix A. The nonhomogeneous H\"{o}lder-Young-Brascamp-Lieb
inequality \label{YBL}}

We have mentioned already in the introduction that the H\"{o}%
lder-Young-Brascamp-Lieb inequality gives the possibility to evaluate the
integrals of the form (\ref{holder}) under conditions on integrability
indices of functions $f_{i}.$

The H\"{o}lder-Young-Brascamp-Lieb inequality was clarified and considerably
generalized recently by Ball \cite{Bal}, Barthe \cite{Barthe}, Carlen, Loss
and Lieb \cite{BLL}, and Bennett, Carbery, Christ and Tao \cite{BCCT1}, \cite%
{BCCT}, the end result being of replacing the linear functionals with
surjective linear operators: $l_{j}({x}):S^{m}\;\rightarrow
S^{n_{j}},j=1,...,k$, with $\cap _{1}^{k}ker(l_{j})=\{0\}$.

Following the remarkable exposition of \cite{BCCT}, \cite{BCCT1}, we give
the formulation of this inequality in the way the most relevant to the
context of the present paper (see Theorem 2.1 of \cite{BCCT}).

Let $H$, $H_{1}$, ..., $H_{m}$ be Hilbert spaces of finite positive
dimensions, each being equipped with the corresponding Lebesgue measure;
functins $f_{j}:H_{j}\rightarrow R$, $j=1,...,m,$ satisfy the integrability
conditions \textit{$f_{j}\in L_{p_{j}}$}, $j=1,...,m.$

Theorem A1 below specifies, in terms of certain linear inequalities on 
\begin{equation*}
z_{j}={\frac{1}{{p_{j}}}},\ j=1,\ldots ,n,
\end{equation*}
the ``power counting polytope'' PCP within which the H\"{o}lder inequality
is valid.

\bigskip

\noindent \textbf{Theorem A1}\textit{\label{t:HYBg} \textbf{(H\"{o}%
lder-Young-Brascamp-Lieb inequality)}. Let $l_{j}(x),j=1,...,k$ be
surjective linear transformationss $l_{j}:H^{m}\;\rightarrow H_{j},$ }$%
j=1,...,m$\textit{. Let $f_{j},\ j=1,\ldots ,k$ be functions $f_{j}\in
L_{p_{j}}(\mu (dx)),\ 1\leq p_{j}\leq \infty $ defined on $H_{j}$, where $%
\mu (dx)$ is Lebesgue measure. }

\textit{Then, the H\"{o}lder-Young-Brascamp-Lieb \ inequality }

\textit{\ 
\begin{equation*}
\biggl|\int_{H}\prod_{j=1}^{m}f_{j}(l_{j}(\mathbf{x}))\mu (d\mathbf{x})%
\biggl|\leq K\prod_{j=1}^{m}\Vert f_{j}\Vert _{p_{j}}\leqno(GH)
\end{equation*}
holds if and only if}

\textit{\medskip \noindent (C1) $\quad \displaystyle{\dim
(H)=\sum_{j}z_{j}\dim (H_{j}),\ }$ }

and

\ 

\textit{\noindent (C2) $\quad \displaystyle{\dim (V)\leq \sum_{j}z_{j}\dim
(l_{j}(V)),\ }$for every subspace ${V\subset H.}$ }

\textit{\noindent Given that (C1) holds, (C2) is equivalent to }

\ 

\textit{(C3) $\quad \displaystyle co{\dim }_{H}{(V)\geq \sum_{j}z_{j}co{\dim 
}_{H_{j}}(l_{j}(V)),\ }$for every subspace ${V\subset H.}$}

\bigskip \textit{Here $\dim (V)$ denotes the dimension of the vector space }$%
V$ \textit{and} $co{\dim }_{H}{(V)}$ \textit{denotes the codimension of a
subspace} \textit{${V\subset H.}$ }

\textit{Note also that any two of conditions (C1), (C2), (C3) imply the
third. }

\bigskip

\textbf{Notes:} 1) The domain of convergence (for fixed $(l_{1},...,l_{k})$)
is called ''power counting polytope'' PCP, cf. the terminology in the
physics literature, where this polytope was already known (at least as
integrability conditions for power functions), in the case $n_{j}=1,\forall
j $. Note that a general explicit form of the facets of PCP when $n_{j}>1$
for some $j$, has not been found yet.

2) Besides the rearrangement techniques of \cite{BL}, this challenging
problem has been also approached recently via "mass transport interpolation"
by \cite{Barthe} and via "heat flow interpolation" by \cite{CLL}.

3) Some related an interesting inequalities and an application to an
analysis of integrals involving cyclic products of kernels can be found in 
\cite{BUZ}.

\bigskip

\subsection*{Acknowledgments}

We thank the referee for constructive remarks and suggestions which helped
to improve the presentation.

\end{document}